\documentclass[11pt]{article}
\usepackage{amssymb,amsmath,graphicx,comment,color,hyperref}
 \usepackage[left=1.25in,right=1.25in,top=1.25in,bottom=1.25in]{geometry}
\newcommand{\matlab}{{\sc matlab}}
\newcommand{\grad}{\nabla}
\newcommand{\R}{\mathbb R}
\newcommand{\C}{\mathbb C}
\newcommand{\E}{\mathbb E}
\newcommand{\bd}{\mathrm{bd}}
\newcommand{\tr}{\mathrm{tr}}
\renewcommand{\Re}{\mathrm{Re}}

\newcommand{\M}{\mathcal M}
\renewcommand{\P}{\mathcal{P}}

\newcommand{\beq}{\begin{equation}}
\newcommand{\eeq}{\end{equation}}
\newcommand{\diag}{\mathrm{diag}}
\renewcommand{\bd}{\mathrm{bd}}
\newcommand{\conv}{\mathrm{conv}}
\newcommand{\twomat}[4]{\left [ \begin{array}{cc} #1 & #2 \\
                        #3 & #4 \end{array}\right ] }
\newcommand{\eps}{\epsilon}
\newcommand{\Zeps}{Z_{\eps}}

\title{Local Minimizers of the Crouzeix Ratio:\\ A Nonsmooth Optimization Case Study}
\author{Michael L. Overton}
\begin{document}
\maketitle

\begin{abstract}
Given a square matrix $A$ and a polynomial $p$, the Crouzeix ratio is the norm of the polynomial on the field of values of $A$
divided by the 2-norm of the matrix $p(A)$. Crouzeix's conjecture states that the globally minimal value 
of the Crouzeix ratio is 0.5, regardless of the matrix order and polynomial degree, and it is known
that 1 is a frequently occurring locally minimal value. Making use of a heavy-tailed distribution
to initialize our optimization computations, we demonstrate for the first time that the Crouzeix ratio has many other locally 
minimal values between 0.5 and 1. Besides showing that the same function
values are repeatedly obtained for many different starting points, we also verify that an approximate nonsmooth
stationarity condition holds at computed candidate local minimizers. We also find that the same locally minimal values are often obtained 
both when optimizing over real matrices and polynomials, and over complex matrices and polynomials.
We argue that minimization of the Crouzeix ratio makes a very interesting nonsmooth optimization case study, illustrating
among other things how effective the BFGS method is for nonsmooth, nonconvex optimization.
Our method for verifying approximate nonsmooth stationarity is based on what may be a novel approach to finding approximate
subgradients of max functions on an interval.
Our extensive computations strongly support Crouzeix's conjecture: in all cases, we find that the smallest locally 
minimal value is 0.5.
\end{abstract}

\section{Introduction}

Let $\P_m$ denote the space
of polynomials with complex coefficients and degree at most $m$, 
let $\M^n$ denote the space of $n \times n$ complex matrices, 
and, for $A\in\M_{n}$, let $W(A)$ denote the field of values (numerical range) of $A$, namely
\beq\label{Wdef}
       W(A) = \{v^*Av:  v\in \C^n, \|v\|_{2} = 1 \}.
\eeq
The field of values is a convex, compact set in the complex plane \cite{HorJoh91}.
For $p\in\P_{m}$ and $A\in\M^n$, the Crouzeix ratio is defined to be
\beq\label{ratio}
      f(p,A) = \frac{\|p\|_{W(A)}} {\|p(A)\|_{2}}
\eeq
where the numerator is
\beq\label{numer}
       \max\left \{|p(z)|: z\in W(A)\right\}
\eeq
and the denominator is the 2-norm of the matrix $p(A)$. Consider the optimization problem
\[
       \min_{p\in\P_{m},A\in\M^{n}}  f(p,A).
\]
This is a nonconvex, nonsmooth optimization problem:  at some pairs $(p,A)$, $f$ is not differentiable.
Crouzeix's famous 2004 conjecture in matrix theory \cite{Cro04} 
postulates that the globally minimal value of 
$f$ is 0.5, regardless of $n$ and $m$. It was established in 2017 \cite{CroPal17} that the globally minimal value
is no less than $\sqrt{2}-1\approx 0.414$.

In our previous work with A.~Greenbaum \cite{GreOve18}, we reported on computational results minimizing $f$
over $n\times n$ real Hessenberg matrices $A$ and real polynomials $p$ with degree at most $n-1$, for $n=3$ through
$n=8$. Using the BFGS method, we ran 100 optimization runs 
for each value of $n$, initialized with the entries of $A$ and the coefficients of $p$ set randomly using
the standard normal distribution. We found that almost all optimization runs generated function values converging either to the conjectured globally
minimal value~0.5, or to the evidently locally minimal value 1. In the former case,
making use of the generalized null space decomposition \cite{GugOveSte15}, we confirmed that the computed final $(p,A)$ always
approximated the conjectured globally minimal ``Crabb matrix'' configurations, to be described below,
with field of values a circular disk. In the latter case, using the Schur factorization, we confirmed that
the final $(p,A)$ always approximated an ``ice-cream-cone'' configuration, again to be described below. 
However, we were not sure whether other locally optimal values of the Crouzeix ratio besides
0.5 and 1 might exist, writing ``\ldots for $n = 7$ and 8, restarting BFGS at and near the final computed pairs [for which $f$ was not
close to 0.5 or 1] did not lead to much improvement, suggesting the possibility
that there are other stationary values of $f$ between 0.5 and 1''. 

In this paper, we report on much more extensive computations minimizing $f$, presenting clear evidence that, in fact, there are many other
stationary values of $f$, specifically locally minimal values of $f$, between 0.5 and 1.
Our new computational results also represent perhaps the strongest evidence yet that Crouzeix's conjecture is true. 
We also think that they present an interesting ``case study'' in nonsmooth optimization.
In particular, they illustrate how effective the BFGS method is at finding locally minimal values even when these occur at nonsmooth
stationary points, despite the fact that the method uses only function and gradient information, not subgradient information. They also illustrate how, 
once candidate stationary points have been identified, approximate stationarity can be verified numerically even at nonsmooth stationary points, 
using approximate subgradient information computed at these points.

The paper is organized as follows. In the next section, we explain how nonsmoothness of the Crouzeix ratio $f$ arises.
In Section  \ref{sec:crabb} we describe the Crabb matrix and ice-cream-cone configurations associated with the stationary
values $0.5$ and $1$  respectively. We briefly describe the computational model in Section \ref{sec:model} and 
explain how to verify approximate
nonsmooth stationarity in Section \ref{sec:stationary}.  We present our experimental results
in Section \ref{sec:results}, and make some concluding remarks in Section \ref{sec:conclude}.

\section{Nonsmoothness of the Crouzeix ratio}

The Clarke subdifferential (generalized gradient) \cite{Cla75} of a locally Lipschitz function
$h$ mapping a Euclidean space to $\R$, evaluated at a point $x$, is the convex hull of the gradient limits
\[
           \partial h(x) = \conv \left\{\lim_{x_{k}\to x} \grad h(x_{k}) \right\},
\] 
where the limit is taken over all sequences $(x_{k})$ converging to $x$ on
which $h$ is differentiable.  Elements of $\partial h(x)$ are called (Clarke) subgradients.
Clearly, if $h$ is continuously differentiable at $x$, then $\partial h(x) = \{\grad h(x)\}$.
If $0\in\partial h(x)$, we say that $x$ is a (Clarke) stationary point of $h$. If, in addition, $h$ is differentiable at $x$ with
$\partial h(x) = \{\grad h(x)\}$, we say that $x$ is a smooth stationary point; otherwise, it is a nonsmooth stationary point.

By identifying $p\in\P_{m}$ with its vector of coefficients $[c_{0},\ldots,c_{m}]^{T}\in\C^{m+1}$, we can view the Crouzeix ratio $f$
given in \eqref{ratio} as a function mapping the Euclidean space $\E=\C^{m+1}\times\M^{n}$,
with real inner product
\[
     \left \langle (c,A),(d,B)\right \rangle = \Re\big ( c^*d + \tr(A^*B)\big),
\]
to $\R$, where$~^*$ denotes complex conjugate transpose. The function $f$ is locally Lipschitz on the set of 
pairs $(p,A)$ for which $p(A)$ is not zero.

By the maximum modulus theorem, the maximum value of $|p(z)|$ on $W(A)$ must be attained on $\bd~W(A)$,
the boundary of the field of values --- and only there, unless $p$ is constant on $W(A)$, which can only occur
if $p$ is constant or $A$ is a multiple of the identity matrix, cases that are of no interest. 
The most important source of nonsmoothness of $f$ is that
\beq \label{Zdef}
            Z(p,A) =\left\{z\in \bd~W(A): \|p\|_{W(A)} = |p(z)|\right\}
\eeq
may contain multiple points. In particular, this occurs at the conjectured global minimizers described in the next section.

As explained in \cite{GreOve18}, there are two other possible sources of nonsmoothness in $f$.
One possibility is that even if $\|p\|_{W(A)}$ is attained only at a single point $z \in \bd~W(A)$, the equation
$z=v^*Av$ in \eqref{Wdef} holds for two or more linearly independent unit vectors $v$.
The other possibility is that the maximum singular value of $p(A)$, which
defines the denominator of the Crouzeix ratio, has multiplicity two or more. 
Since neither of these cases occurs either at the conjectured global minimizers or at the apparent local minimizers found
in our computations, we will not consider them further.

\section{The Crabb matrix, conjectured global minimizers, and ice-cream-cone stationary points}\label{sec:crabb}

Pairs $(\tilde p,\tilde A)$ for which the Crouzeix ratio is 0.5 are known. Given an integer $k$ with $2\leq k\leq \min(n,m+1)$,
define the polynomial $\tilde p \in \P_m$ by $\tilde p(\zeta)=\zeta^{k-1}$,
set the matrix $\Xi_k\in\M^k$  to
\beq \label{choidef}
   \twomat{0}{2}{0}{0} \mathrm{~if~}k=2,  \mathrm{~or~} 
     \left [ \begin{array}{ccccccc}
                 0 & \sqrt{2} &        &       &       &      &           \\
                         & \cdot    & 1      &       &       &      &           \\
                         &          & \cdot  & \cdot &       &      &           \\
                         &          &        & \cdot & \cdot &      &           \\
                         &          &        &       & \cdot & 1    &           \\
                         &          &        &       &       & \cdot & \sqrt{2} \\
                         &          &        &       &       &       &  0 \end{array} \right ]
\mathrm{~if~}k>2,
\eeq
and set $\tilde A = \diag(\Xi_k, 0) \in \M^n$. The matrix $\Xi_{k}$ was called the Choi-Crouzeix matrix of order $k$ in \cite{GreOve18},
but after the paper was published, A.~Salemi\footnote{Private communication, 2017}
 informed us that it was introduced much earlier in a different context by Crabb \cite{Cra71}.
The field of values of the Crabb matrix $\Xi_{k}$ is the unit disk, so the numerator of the Crouzeix ratio is 1, and
$\tilde p(\tilde A)=\tilde A^{k - 1}=\diag(\Xi_k^{k-1},0)$ is a matrix with just one nonzero, namely a 2 in the $(1,k)$ position,
so the denominator is 2; hence, the ratio is 0.5.

Since $|\tilde p|$ is constant on the boundary of the unit disk, we have that 
\[
           Z(\tilde p,\tilde A) = \{z\in\C: |z|=1\},
\]
the unit circle, resulting in nonsmoothness of the Crouzeix ratio $f$ at $(\tilde p,\tilde A)$.
Together with A.~Greenbaum and A.S.~Lewis \cite{GreLewOve17}, we derived the Clarke subdifferential of $f$ at $(\tilde p,\tilde A)$, and
established that $(\tilde p,\tilde A)$ is a nonsmooth stationary point of $f$. The analysis also showed that 
$f$ is directionally differentiable at $(\tilde p,\tilde A)$ and that the directional derivative of $f$ is nonnegative in every direction in $\E$.
Although this does not imply that $(\tilde p,\tilde A)$ is even a local minimizer of $f$, a
discussion of how one might extend the analysis towards the goal of proving local (but not global!) minimality is given
in \cite[p.\ 242]{GreLewOve17}.

The property that $f(\tilde p, \tilde A)=0.5$ extends easily to pairs $(p,A)$ where $p(\zeta)=(\zeta-\lambda)^{k-1}$ and
\[
                 A = \lambda I + \beta \, U\diag(\Xi_k,B)U^*,
\]
for any $\lambda\in\C$, nonzero $\beta\in\C$,  unitary matrix $U$, and matrix $B$ with $W(B)$ contained in the unit disk. 
We conjecture that such pairs $(p,A)$ are the \emph{only} ones for which $f(p,A)=0.5$.
However, if the condition that $p$ is a polynomial is relaxed to allow it to be any analytic function, there are
many choices for $(p, A)$ for which the ratio is 0.5; for the case $n=3$, see \cite[Sec.\ 10]{Cro16}.

In our computations, we often find the locally minimal value 1. This occurs when
the matrix $A$ is block diagonal of the form $A=\diag(\lambda, B)$, with $\lambda\in\C$, 
$\lambda\not\in W(B)$. In this configuration, $W(A)=\conv(\lambda, W(B))$ with $\bd~W(A)$ consisting only of
$\lambda$, part of $\bd~W(B)$ and two line segments connecting $\lambda$ to $W(B)$. Hence, $W(A)$ has a vertex at $\lambda$, and 
often has the appearance of an ice cream cone, as illustrated by the examples reported in \cite[Fig.~4]{GreOve18}. 
If, in addition,
\[
                |p(\lambda)|   > |p(\nu)| \mathrm{~for~all~}\nu\in W(A), \nu\not = \lambda,  
\]
so that that $Z(p,A)$ consists only of the single point $\lambda$, and
\[
               |p(\lambda)|  > \|p(B)\|_{2},
\]
it is immediate that both the numerator and denominator of the Crouzeix ratio are $|p(\lambda)|$, so $f(p,A)=1$.
Furthermore, we showed in \cite[Thm.~2]{GreOve18} that $f$ is differentiable 
at such $(p,A)$ and that its gradient is zero, so that $(p,A)$ is a smooth
stationary point (though not necessarily a local minimizer).

\section{The computational model}\label{sec:model}

We use the same computational model as in \cite{GreOve18}, so we briefly explain only the main points.
It is well known  \cite{Kip51} that $\bd~W(A)$, the boundary of $W(A)$, can be
characterized as
\beq\label{bdWchar}
      \bd~W(A) = \left \{z_\theta = v_\theta^*A v_\theta : \theta \in [0,2\pi) \right \}
\eeq
where $v_\theta$ is a normalized eigenvector corresponding to the largest eigenvalue of the
Hermitian matrix
$$
      H_\theta = \frac{1}{2}\left (e^{i\theta} A + e^{-i\theta} A^*\right ).
$$
The proof uses a supporting hyperplane argument \cite[Prop.\ 2]{GreLewOve17}.
To accurately and efficiently approximate $\bd~W(A)$, we use Chebfun \cite{Chebfun14},
a system for approximating functions on a real interval to machine precision accuracy by adaptive
Chebshev approximation. Chebfun's function {\tt fov} computes a complex-valued ``chebfun'' approximating the extreme points of
$\bd~W(A)$ on the interval $[0,2\pi]$, generating interpolation points $\theta\in (0,2\pi)$  automatically. 
A second output argument returns any line segments in the boundary as well.\footnote{This modification to {\tt fov} was written by the author.}


For the optimization calculations, we use the BFGS method,
devised independently in 1970 by Broyden, Fletcher, Goldfarb and Shanno for unconstrained optimization of
differentiable functions, but which is also extremely effective for nonsmooth optimization \cite{LewOveNSOquasi}.
BFGS requires the computation of the Crouzeix ratio $f(p,A)$ and its gradient at a sequence
of iterates generated by the method. The main cost in computing $f(p,A)$ is that of constructing the chebfun representing 
$\bd~W(A)$.  Computing the numerator of \eqref{ratio} is then
done by invoking two \matlab\ functions that have been overloaded
to be applicable to chebfuns, namely {\tt polyval} and {\tt norm(.,inf)}, while computation of the denominator, the 2-norm of
$p(A)$, is carried out by calls to two standard \matlab\ functions, {\tt polyvalm} and {\tt norm}.
Once $f(p,A)$ has been computed, the additional computation required to obtain its gradient is minimal,
even though the formula is complicated: see \cite[Theorem~1]{GreOve18} for details.
In order to compute the gradient, we need to know $Z(p,A)$, which tells us where $\|p\|_{W(A)}$ is attained;
 this information is returned by Chebfub's {\tt norm(.,inf)} function.
A natural question is: what is the method to do if $Z(p,A)$ contains multiple points, and hence $f$ is not differentiable at $(p,A)$? 
The answer is that since BFGS uses a ``gradient paradigm'', not a ``subgradient paradigm'' \cite{LewOveNSOquasi,AslOve21}, this
possibility, which is essentially impossible to check exactly in finite precision, is simply ignored.
In practice, the algorithm will virtually never compute pairs $(p,A)$ where $Z(p,A)$ contains multiple points, except in the limit. 
Clearly, small changes in $(p,A)$ may result in large changes in the computed gradient, but this is 
inherent in nonsmooth optimization, and in fact explains to a large extent why BFGS works so well 
in this context \cite[p.~130]{LewOveNSOquasi}.
The BFGS method is a line-search descent method, meaning that at every iteration it
uses an inexact line search to repeatedly evaluate the minimization objective $f$ along a search direction in the variable space until the
so-called Armijo-Wolfe conditions are satisfied. If, due to rounding errors, it is not possible to satisfy these conditions in a reasonable
number of steps, the BFGS method is terminated.

\section{Verifying approximate nonsmooth stationarity} \label{sec:stationary}

In order to verify approximate nonsmooth stationarity of computed pairs $(p,A)$ \emph{a posteriori}, 
we need to consider some sort of approximate subdifferential of $f$ at $(p,A)$. 
It is well known \cite[Thm.~10.31]{RocWet98} that the subdifferential of a max function on
an interval is the convex hull of gradients of the component functions evaluated at points where the max is attained. Hence the importance
of the set $Z(p,A)$, which underlies the derivation of the subgradients of the Crouzeix ratio at the Crabb matrix configuration 
$(\tilde p,\tilde A)$ given in \cite{GreLewOve17}.
We now introduce what may be a novel idea for approximating the subdifferential of a max function on an interval: instead of
the convex hull of the gradients evaluated where the max is attained exactly, we use the convex hull of gradients evaluated
at \emph{local} maximizers for which the locally maximal value is sufficiently close to the globally maximal value.
For $\eps \geq 0$, define
\begin{align} \label{localmax}
           \Zeps (p,A) =  \{z \in \bd~W(A): &  ~z \text{ is a local maximizer of } \zeta \mapsto |p(\zeta)| \\
                                                         &  \text{ and } |p(z)| \geq (1-\eps) \|p\|_{W(A)} \}.\nonumber
\end{align}
Clearly, $Z_{0}(p,A)=Z(p,A)$.  Then, replace $Z(p,A)$ in \cite[Eq.\ (9)]{GreLewOve17} by $\Zeps(p,A)$ and use the resulting modified formula
for $\partial f(p,A)$ in \cite[Thm.\ 3]{GreLewOve17} as our approximate subdifferential, say $\partial^{\,\eps}f(p,A)$. 
This should not be confused
with other usages of the word ``approximate'' in subdifferential analysis which sometimes mean ``limiting'' \cite[p.\ 347]{RocWet98},
or by approximating the subdifferential using perturbations to $(p,A)$ \cite{Gol77,BurLewOveMOR,BurLewOveGradSamp}. 

Given a computed pair $(p,A)$, we use Chebfun's {\tt max(.,'local')} function to compute all local
maximizers of $|p|$ on $\bd~W(A)$. If $|p|$ is constant on $\bd~W(A)$, which can only happen if $W(A)$ is a disk, as in the
case that $A$ is a Crabb matrix, Chebfun returns no local maximizers, so we
forgo computing the stationarity measure in this case. At all other local minimizers of $f$ such
as those described below, the number of local maximizers of $|p|$ on $\bd~W(A)$ is necessarily finite, and hence the set $\Zeps(p,A)$ is
a discrete set, as opposed to a continuum that would be obtained if we eliminate the ``local maximizer'' condition in
\eqref{localmax}. Note also that if $W(A)$ has
the ice-cream-cone configuration, with $|p|$ maximized on $\bd~W(A)$ only at the vertex $\lambda$, then the local maximizers
returned by Chebfun do not include $\lambda$, since $\bd~W(A)$ is not smooth there. Hence, it is important to compute the
global maximum using Chebfun's {\tt max(.)} as well as the local maximizers with {\tt max (.,'local')}.

Finally, since the exact nonsmooth stationarity condition is $0\in\partial f(p,A)$, we compute
\beq \label{ddef}
           d = \mathrm{argmin} \{\|g\|_{2} : g\in\partial^{\,\eps}f(p,A)\},
\eeq
the solution of a convex quadratic programming problem, and use $\|d\|_{2}$ as a measure of approximate nonsmooth stationarity.

\section{The new computational results}\label{sec:results}

\begin{figure}
\begin{center}
\includegraphics[scale=0.9,trim={1.5cm 1.5cm 1.4cm 0.5cm},clip]{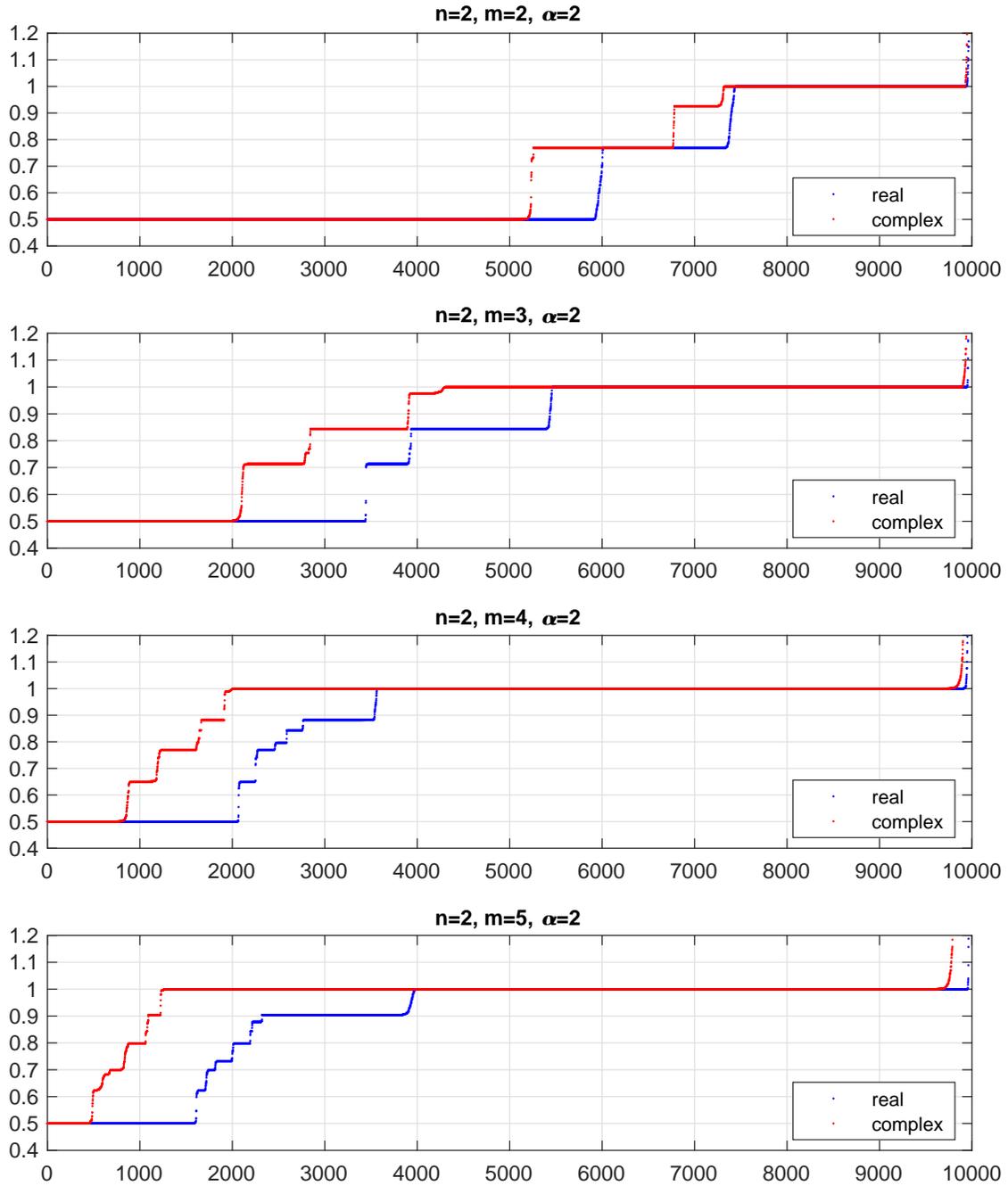}
\end{center}
\caption{
Sorted final values of the Crouzeix ratio $f$ obtained from 10,000 runs for $n=2$, $m=2,3,4,5$,
optimizing over \emph{real} $n\times n$ matrices and \emph{real} polynomials of degree \emph{at most} $m$ (blue dots)
\emph{and} optimizing over \emph{complex} $n\times n$ matrices and \emph{complex} polynomials of degree \emph{at most} 
$m$ (red dots), with starting points generated by the heavy-tail distribution defined by \eqref{heavy}, with $\alpha=2$.
Note that many of the blue dots are overwritten by red dots.}
\label{fig:n2m2345}
\end{figure}

M.~Hairer\footnote{Private communication, 2019} pointed out that by starting with initial data generated from the normal distribution, as done in the results reported in \cite{GreOve18},
we might be biasing the optimization towards matrices whose fields of values are disks. He suggested trying initial data generated from 
distributions with heavy tails, e.g., numbers of the form
\begin{equation}\label{heavy}
           x e^{\alpha x^{2}}
\end{equation}
where $x$ is obtained from the standard normal distribution and $\alpha>0$. Indeed, experiments show that 
for fixed $n$, as $\alpha$ increases, it becomes much less likely that BFGS will generate Crabb matrix configurations with fields
of values a disk and Crouzeix ratio $0.5$. But instead of finding lower values that would disprove the conjecture, 
what happens instead is that BFGS is much more likely to find other locally minimal values of the Crouzeix ratio between
0.5 and 1. In the experiments reported here, we fix $\alpha=2$ in \eqref{heavy}. This is large enough to 
discover many more locally minimal values than using $\alpha = 0$ or $\alpha = 1$, but not so large that we have problems with overflow.

In the results reported here, in addition to optimizing over real polynomials and real matrices, we also present results obtained by optimizing
over complex polynomials and complex matrices. Without loss of generality, since the Crouzeix ratio is invariant to unitary similarity 
transformations of the matrix, in the real case we restrict the matrices to upper Hessenberg
form; in the complex case, we restrict them to upper triangular form. Furthermore, in the real case, since the field of values is symmetric
with respect to the real axis, we compute its boundary only in the closed upper half-plane, restricting $\theta$ in \eqref{bdWchar}
to the interval $[\pi,2\pi]$, and modifying the definitions of $Z(p,A)$ and $\Zeps(p,A)$ to restrict them to points $z$ with $\Re(z) \geq 0$.

\subsection{The case $n=2$}

Figure~\ref{fig:n2m2345} shows the final values of the Crouzeix ratio $f$ for $n=2$ with 
the maximal degree $m$ ranging over 2, 3, 4 and 5, obtained by running 
BFGS\footnote{Using {\sc hanso} 2.2 (www.cs.nyu.edu/overton/hanso), with {\tt options.normtol = 1e-8}.}
initialized from 10,000 randomly 
generated starting points using \eqref{heavy} with $\alpha=2$, for each of the real and complex cases. 
The blue dots show the \emph{sorted} final values of the Crouzeix ratio obtained
from optimizing over real matrices of order $n$ and real polynomials with maximal degree $m$, and the red dots show 
the \emph{sorted} final values when optimizing over complex
matrices of order $n$ and complex polynomials with maximal degree $m$. 
Note that since the results for the real case are plotted first, many of the blue dots are overwritten by red dots.
The plateaus clearly indicate locally minimal values, as  these values are found repeatedly from
many starting points. Furthermore, we see that most of the locally minimal values found are the \emph{same} for both the real and
complex cases. That said, there is not a great deal of similarity between the real and complex results. The width
of the plateaus of locally minimal values found represents the frequency with which they are found, and this varies greatly
between the real and complex cases. Partly for this reason, it is difficult to check whether there is
a one-to-one correspondence between the locally minimal values found in the real and complex cases;
there are some possible counterexamples if we assume that 10,000 runs is enough to see these
features, but obviously we have no basis for this. The locally minimal values 0.5 and~1 are clearly apparent in every case.
The former, 0.5, is known to be globally minimal for $n=2$, and is globally minimal for all $n$ and $m$ if Crouzeix's conjecture
is true. The latter value, 1, is the stationary value associated with ice-cream-cone configurations discussed above.

\begin{table}
\begin{center}
\small
\begin{tabular}{|c|c|c|c|c|c|c|}
\hline\hline
   & run \# & numer & denom & $f$ & $|Z_\epsilon|$ & $\|d\|$ \\ 
\hline\hline 
R & 1 & 4.195e+03 & 8.391e+03 &  0.5000000000 &  & \\ 
R & 3400 & 1.031e+05 & 2.061e+05 &  0.5000188514 &  & \\ 
\hline 
R & 3500 & 7.808e+06 & 1.095e+07 &  0.7132185867 & 1 & 3.313e-08 \\ 
R & 3850 & 2.371e+11 & 3.325e+11 &  0.7132189899 & 1 & 4.689e-09 \\ 
\hline 
R & 4000 & 7.904e+11 & 9.368e+11 &  0.8437496323 & 1 & 8.131e-09 \\ 
R & 5200 & 9.164e+17 & 1.086e+18 &  0.8437501759 & 1 & 5.790e-09 \\ 
\hline 
R & 5500 & 5.624e+00 & 5.624e+00 &  1.0000000000 & 1 & 8.436e-16 \\ 
R & 9800 & 6.694e+14 & 6.694e+14 &  1.0000000084 & 1 & 4.322e-09 \\ 
\hline \hline
C & 1 & 1.362e+15 & 2.724e+15 &  0.5000000000 & & \\ 
C & 1900 & 5.001e+03 & 9.996e+03 &  0.5002960578 & & \\ 
\hline 
C & 2200 & 3.113e+07 & 4.365e+07 &  0.7132064490 & 2 & 1.574e-05 \\ 
C & 2700 & 1.122e+10 & 1.574e+10 &  0.7132194699 & 2 & 1.263e-07 \\ 
\hline 
C & 2900 & 7.691e+08 & 9.115e+08 &  0.8437482418 & 2 & 1.809e-07 \\ 
C & 3800 & 2.371e+11 & 2.810e+11 &  0.8437500721 & 2 & 5.533e-09 \\ 
\hline
C & 4500 & 2.470e+07 & 2.470e+07 &  1.0000000000 & 3 & 1.845e-04 \\ 
C & 9800 & 7.825e+18 & 7.825e+18 &  1.0000004701 & 1 & 3.734e-09 \\ 
\hline
\end{tabular}
\end{center}
\caption{Four locally minimal values for $n = 2$, $m = 3$, including 0.5 and 1. 
The first column indicates whether the data is from the real or complex run, and the second column shows the relevant run number. 
The next columns show the numerator and denominator of the final Crouzeix ratio, as well as the ratio $f$ itself.
The final two columns show the number of points in $\Zeps(p,A)$ and the resulting
approximate stationarity measure, using $\eps=10^{-4}$.
}
\label{tab:n2m3}
\end{table}

It's interesting to see how the other locally minimal values visible in Figure~\ref{fig:n2m2345} vary with $m$. 
In both the real and complex cases, in the case $m=3$,
the widest plateau between 0.5 and 1  has value 0.84375, and this is also visible as shorter
plateaus in the cases $m=4$ and $m=5$, but it does not appear in the case $m=2$.
Meanwhile, again in both the real and complex cases, in the case $m=2$, the widest plateau between 0.5 and 1, with value
0.7698, is reduced to much narrower plateaus in the case $m=3$, $m=4$ and $m=5$, which are only visible with \matlab's
zoom tool. This behavior is unexpectedly complicated.

Table~\ref{tab:n2m3} shows more details associated with four locally minimal values clearly visible as plateaus in the 
second panel of Figure~\ref{fig:n2m2345}, for $n=2$, $m=3$, namely, 0.5, 0.713, 0.844 and~1. 
These locally minimal values are found in both the real case (indicated by R in the first column of the table)
and the complex case (indicated by C in the table).\footnote{The complex
runs also identify a fifth locally minimal value, 0.977, but it is not clear whether this is a locally minimal
value in the real case; in any case, we cannot conclude that from the figure.}
The second column of the table indicates the run number in the horizontal scale used
in Figure~\ref{fig:n2m2345}. These indices are chosen to \emph{roughly} correspond to the first and last points in each plateau.
Consequently, corresponding pairs of values of the Crouzeix ratio $f$, shown in the fifth column,
indicate the approximate precision to which the locally minimal values are found.
The third and fourth columns show the associated numerator and denominator of
\eqref{ratio}. Note that, despite the enormous values of the numerator and denominator
of the Crouzeix ratio, the ratio $f$ is consistently computed to several digits of agreement even over
large numbers of starting points. For example, the first 3400 real runs approximate 0.5 to 4 digits of accuracy, while
the first 1900 complex runs approximate 0.5 to 3 digits.
The final two columns of Table~\ref{tab:n2m3} show two quantities associated with the
approximate stationarity measure described in Section \ref{sec:stationary}, namely, the number of points in $\Zeps(p,A)$ and 
the  2-norm of the corresponding vector obtained in \eqref{ddef}, using $\eps=10^{-4}$.
For the runs which find $f$ approximately equal to 0.5, these values are omitted,
because, at the Crabb matrix configurations, $Z(p,A)$ is a continuum, and, as noted earlier, Chebfun does not find
any local maxima of $|p|$ on $\bd~W(A)$ in this case. In all other cases, the value of $\|d\|$ is small, strongly indicating
that the computed $(p,A)$ is an approximate stationary point. We comment further below on the number of points in $\Zeps(p,A)$.

\begin{figure} 
\begin{center}
\includegraphics[scale=0.75,trim={4.5cm 3cm 3cm 1.5cm},clip]{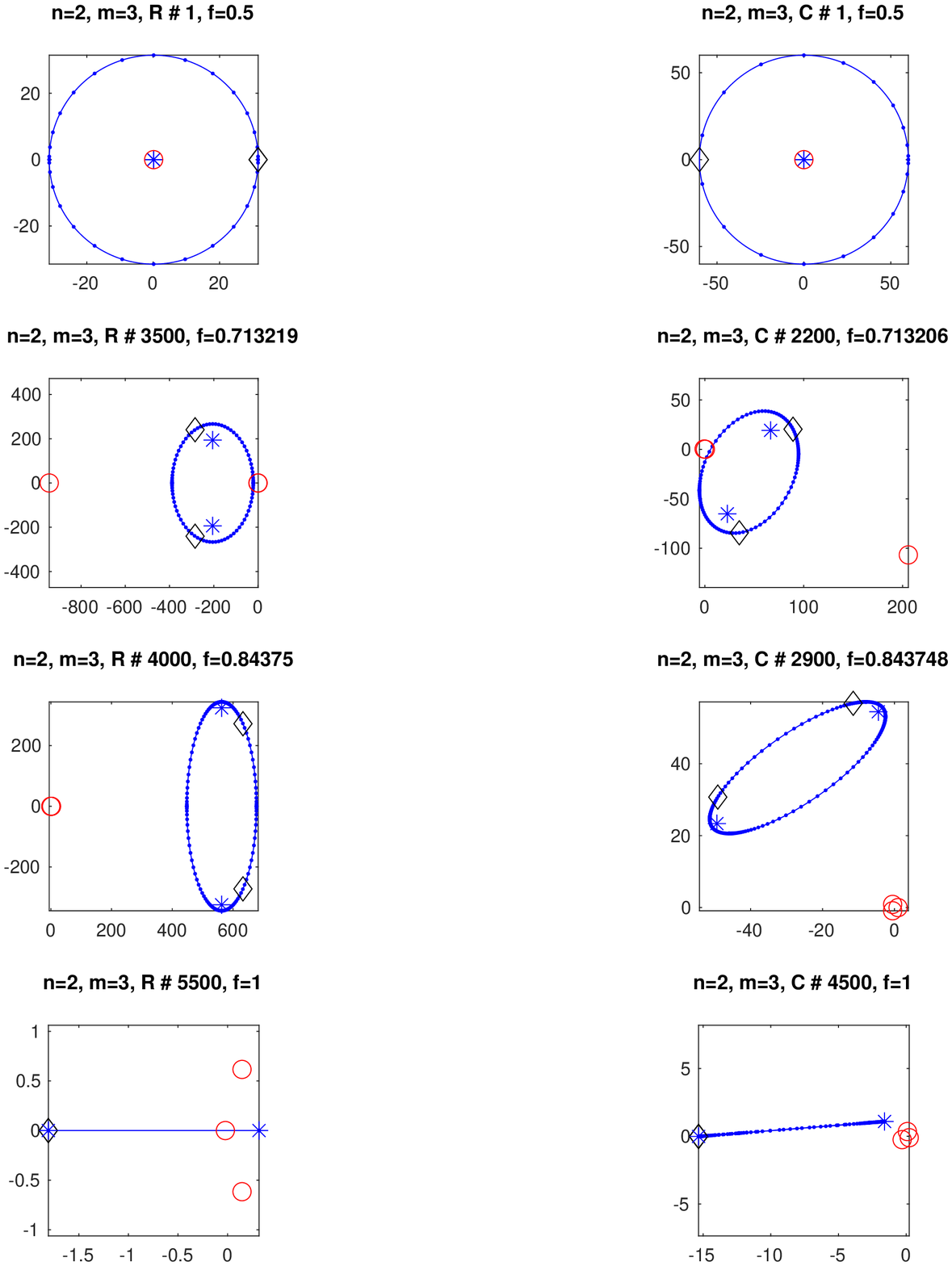}
\end{center}
\caption{Configurations of some local minimizers obtained for $n=2$ and $m=3$ with locally minimal
values 0.5, 0.713, 0.844 and 1, corresponding to real runs 1, 3500, 4000 and 5500 on the left and 
complex runs 1, 2200, 2900 and 4500 on the right.
Solid blue curves show the boundaries of the field of values $W(A)$,
blue asterisks show the eigenvalues of $A$, 
red circles show the roots of $p$, and black diamonds show the points on $W(A)$ where $\|p\|_{W(A)}$ is 
attained.}
\label{fig:n2m3fovs-subplot}
\end{figure}

Figure~\ref{fig:n2m3fovs-subplot} shows the fields of values of $A$ for the final configuration obtained
by the first real run (left)  and first complex run (right) given in Table~\ref{tab:n2m3}
associated with the locally minimal values 0.5, 0.713, 0.844  and 1
(top row, second row, third row and bottom, respectively). In the real case, the fields of values are necessarily symmetric
w.r.t.\ the real axis. Blue asterisks denote the eigenvalues of $A$, red circles show the roots of $p$, and
black diamonds show the points in $\Zeps(p,A)$, where $\|p\|_{W(A)}$ is (nearly) attained.
Clearly, when shifted, scaled and rotated, the final fields of values found by the complex runs are 
very similar to the ones found by the real runs that approximate the same locally minimal value. 

The known configuration attaining $f=0.5$ is $(p,A)$ where $A$ is a Jordan block (the Crabb matrix for $n=2$),
for which $W(A)$ is a disk, and $p(z)=z-\lambda$, where $\lambda$ is the double eigenvalue of $A$. 
We see this configuration in the top two panels of Figure~\ref{fig:n2m3fovs-subplot}.  Note that the two eigenvalues 
of the final computed matrix $A$ and one of the roots
of the final computed polynomial $p$ (nearly) coincide at the center of the disk.
Since the computations are being done with maximal degree $m=3$, there are two other roots of $p$, which,
for the optimal $p$, are $\infty$, but, in our computations, are numbers with very large modulus that are not shown. 
In theory, the optimal $\|p\|_{W(A)}$ is attained at every point on the boundary of $W(A)$, but
only one point is shown.

In the second row of Figure~\ref{fig:n2m3fovs-subplot}, the eigenvalues of $A$ are distinct, so $W(A)$ is elliptical,
but two of the three roots of $p$ (nearly) coincide, a little outside $W(A)$, while the third root is on the other side of 
and further away from $W(A)$. In the third row, $W(A)$ is a more
eccentric ellipse than in the second row, and the three roots of $p$ are all (nearly) coincident.\footnote{Since
optimizing over complex matrices and complex polynomials requires much more computation than in the real case, the
results are likely less accurate in the complex case, as is suggested by the three roots of $p$ being less close to coincident.}
These are interesting, and decidedly non-random, configurations.
In the second and third rows, $\|p\|_{W(A)}$ is attained at two points on the boundary of $W(A)$. 
In the real cases, this is a consequence of
the imposed structure: since $\|p\|_{W(A)}$ is attained at a complex point, it must also be attained at the conjugate point.
For this reason, as shown in Table~\ref{tab:n2m3}, $\Zeps(p,A)$ contains only one point, since only points in the closed upper half-plane 
are admissible in the real case. It follows that the minimizer is a smooth stationary point in the real data space, and the
vector $d$ whose norm is shown in the table is actually the gradient. 
On the other hand, in the complex case, no such structure is imposed,
and the double attainment is reflected by $\Zeps(p,A)$ having two points in this case. 
Hence, this is a nonsmooth minimizer in the complex data space, and the vector $d$ whose norm is shown in Table~\ref{tab:n2m3}
is an approximate subgradient, not a gradient.

The final configuration for the locally optimal value 1 is a diagonal matrix whose
field of values is a line segment, which, as $n=2$, is a degenerate case of the ice-cream-cone fields of values 
mentioned earlier. In the bottom row of Figure~\ref{fig:n2m3fovs-subplot}, in the real case, indeed we see that $W(A)$ is a line segment,
so $\Zeps(p,A)$ has a single point, $f$ is differentiable and the vector $d$ is a gradient, with small norm.
On the other hand, although the field of values computed in the complex case is nearly a line segment, it is not exactly a
line segment, and $\Zeps(p,A)$ actually has 3 points, although they are all nearly identical. Consequently, the vector $d$ is
an approximate subgradient, but it still has small norm, showing the robustness of the calculations.

\begin{figure}
\begin{center}
\includegraphics[scale=0.9,trim={1.5cm 1.5cm 1.4cm 0.5cm},clip]{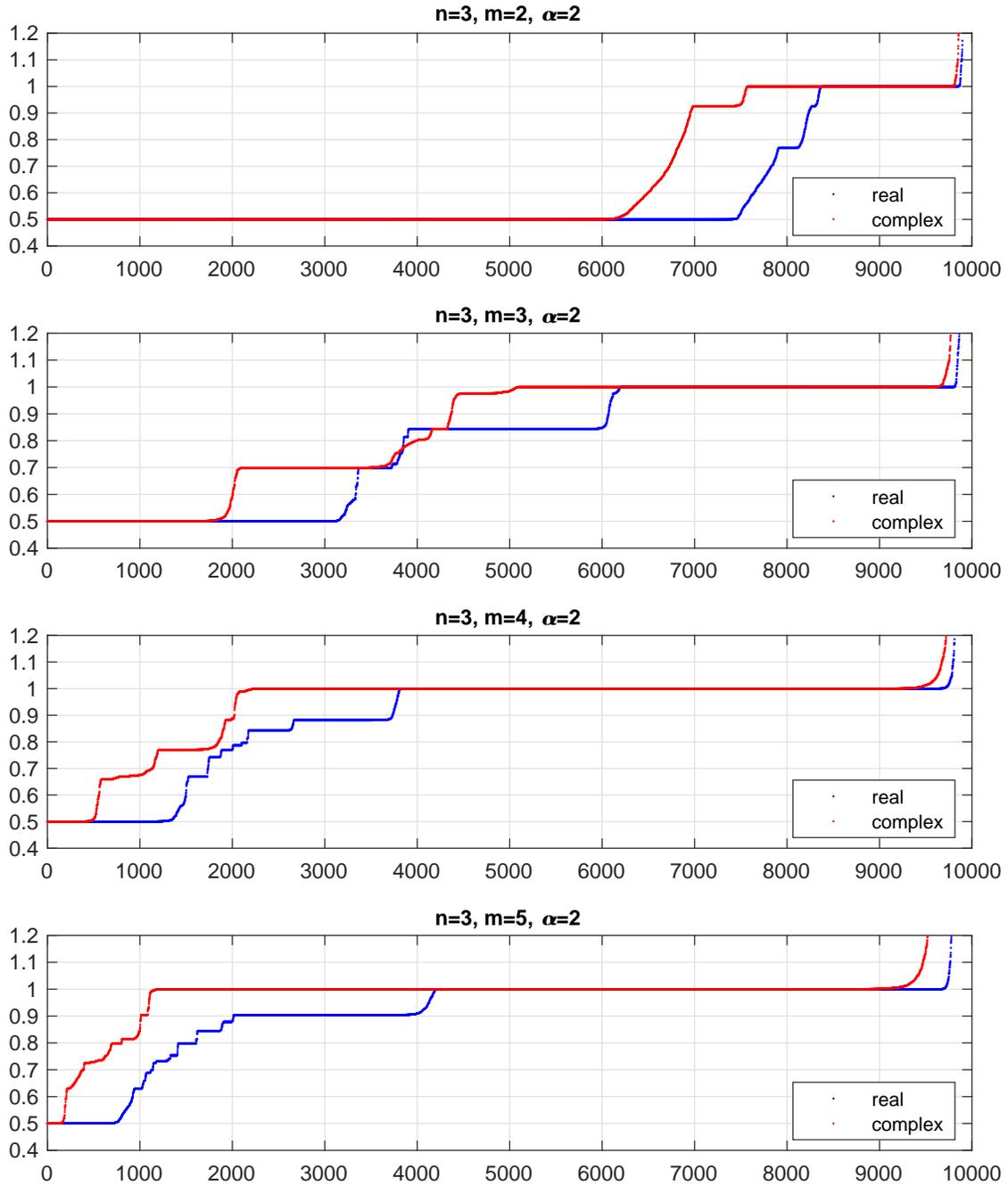} 
\end{center}
\caption{
Sorted final values of the Crouzeix ratio $f$ obtained from 10,000 runs for $n=3$, $m=2,3,4,5$,
optimizing over \emph{real} $n\times n$ matrices and \emph{real} polynomials of degree \emph{at most} $m$ (blue dots)
\emph{and} optimizing over \emph{complex} $n\times n$ matrices and \emph{complex} polynomials of degree \emph{at most} 
$m$ (red dots), with starting points generated by the heavy-tail distribution defined by \eqref{heavy}, with $\alpha=2$.
Note that many of the blue dots are overwritten by red dots.}
\label{fig:n3m2345}
\end{figure}

 \begin{table}
\begin{center}
\small
\begin{tabular}{|c|c|c|c|c|c|c|}
\hline\hline
    & run \# & numer & denom & $f$ & $|Z_\epsilon|$ & $\|d\|$ \\ 
\hline\hline 
R & 1 & 1.607e+14 & 3.213e+14 &  0.5000000000  & & \\ 
R & 3000 & 4.982e+02 & 9.962e+02 &  0.5001593709 & & \\ 
\hline 
R & 3400 & 5.766e+09 & 8.264e+09 &  0.6978015654 & 2 & 6.937e-04 \\ 
R & 3700 & 3.127e+17 & 4.481e+17 &  0.6978024061 & 2 & 1.022e-05 \\ 
\hline 
R & 4000 & 2.395e+12 & 2.839e+12 &  0.8437498418 & 1 & 1.370e-08 \\ 
R & 5800 & 6.651e+16 & 7.883e+16 &  0.8437562126 & 2 & 3.681e-08 \\ 
\hline 
R & 6250 & 2.489e+03 & 2.489e+03 &  1.0000000000 & 1 & 6.227e-12 \\ 
R & 9800 & 7.920e+19 & 7.920e+19 &  1.0000670371 & 1 & 3.732e-09 \\ 
\hline \hline
C & 1 & 4.262e+13 & 8.523e+13 &  0.5000000000 & & \\ 
C & 1400 & 5.786e+10 & 1.157e+11 &  0.5001167444 & & \\ 
\hline 
C & 2200 & 7.823e+10 & 1.121e+11 &  0.6978004851 & 3 & 6.705e-05 \\ 
C & 3200 & 4.851e+19 & 6.950e+19 &  0.6979798248 & 3 & 9.135e-07 \\ 
\hline 
C & 4200 & 4.439e+16 & 5.261e+16 &  0.8437493557 & 2 & 1.055e-08 \\ 
C & 4300 & 6.078e+19 & 7.204e+19 &  0.8437505791 & 2 & 7.263e-09 \\ 
\hline 
C & 5500 & 1.005e+14 & 1.005e+14 &  1.0000000000 & 1 & 3.356e-13 \\ 
C & 9500 & 3.181e+19 & 3.181e+19 &  1.0000027774 & 1 & 5.767e-09 \\ 
\hline
\end{tabular}
\end{center}
\caption{Four locally minimal values for $n = 3$, $m = 3$, including 0.5 and 1. 
The first column indicates whether the data is from the real or complex run, and the second column shows the relevant run number. 
The next columns show the numerator and denominator of the final Crouzeix ratio, as well as the ratio $f$ itself.
The final two columns show the number of points in $\Zeps(p,A)$ and the resulting
approximate stationarity measure, using $\eps=10^{-4}$.
}
\label{tab:n3m3}
\end{table}

\subsection{The case $n=3$}

\begin{figure} 
\begin{center}
\includegraphics[scale=0.75,trim={4.5cm 3cm 3cm 1.5cm},clip]{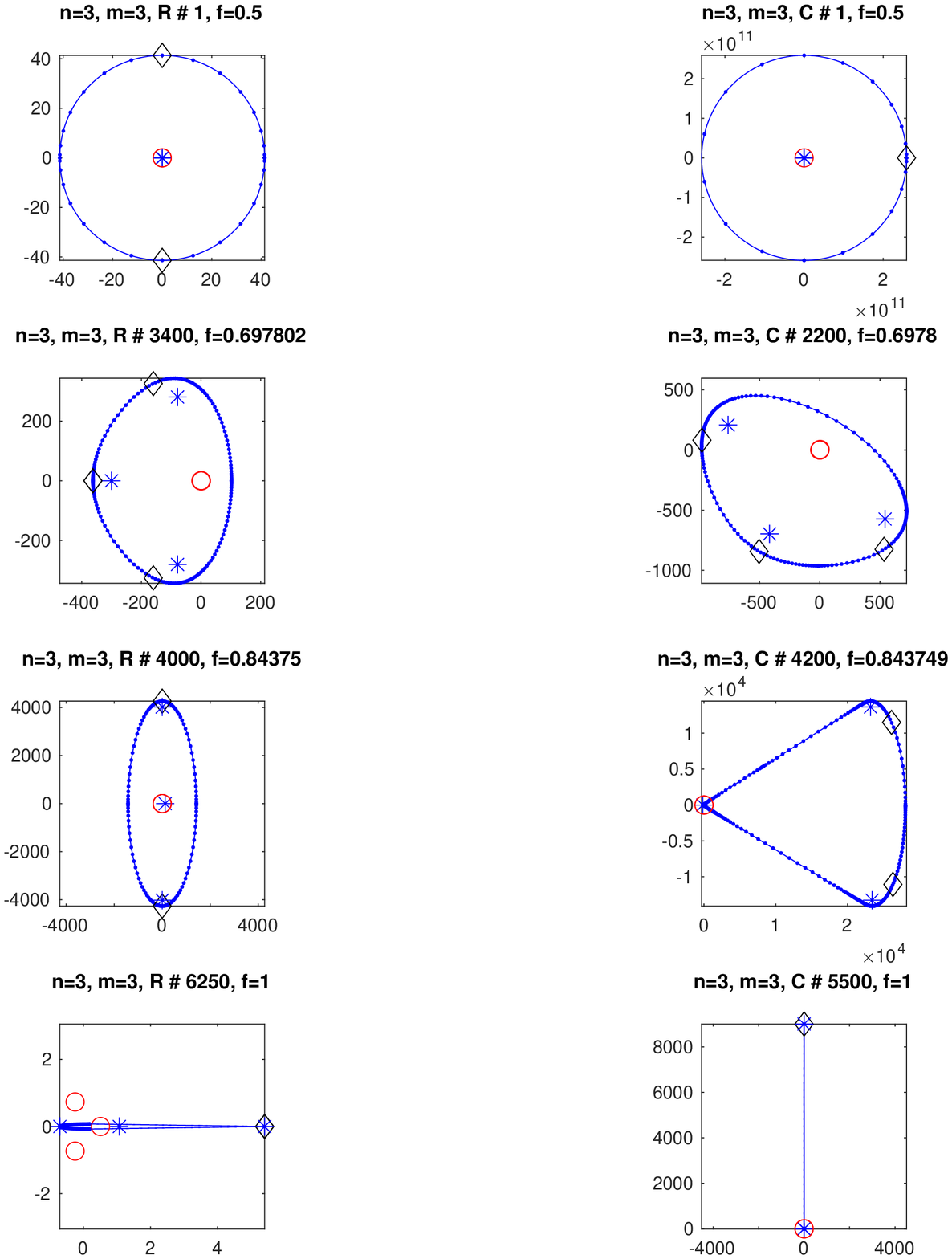}
\end{center}
\caption{Configurations of some local minimizers obtained for $n=3$ and $m=3$ with locally minimal
values 0.5, 0.698, 0.844 and 1, corresponding to real runs 1, 3400, 4000 and 6250 on the left and 
complex runs 1, 2200, 4200 and 5500 on the right.
Solid blue curves show the boundaries of the field of values $W(A)$,
blue asterisks show the eigenvalues of $A$, 
red circles show the roots of $p$, and black diamonds show the points on $W(A)$ where $\|p\|_{W(A)}$ is 
attained.}
\label{fig:n3m3fovs-subplot}
\end{figure}

Figure~\ref{fig:n3m2345} shows sorted final values for the case $n=3$, for $m=2,3,4,5$, again using 10,000 starting points
for each of the real and complex runs. 
Let us again focus on the results for $m=3$, shown in the second panel.
Table~\ref{tab:n3m3} shows details associated with the four locally minimal values 0.5, $0.698$, $0.844$ and 1
that are clearly visible in the figure.\footnote{Again, we omit a fifth locally minimal value, 0.977, that is clearly identified by the
complex runs; unlike in the case $n=2$, $m=3$, zooming in indicates that there is also a small plateau with this value 
for the real runs.}  The locally minimal value $0.698$ is
significantly less than the value $0.713$ observed in Figure~\ref{fig:n2m2345} for the case $n=2$, $m=3$,  
but the locally minimal value $0.844$ is the same as the value observed in the case $n=2$, $m=3$. 
Figure~\ref{fig:n3m3fovs-subplot} shows the fields of values of $A$ for the final configuration obtained
by the first real run (left)  and first complex run (right) given in Table~\ref{tab:n3m3}
associated with the four locally minimal values 0.5, $0.698$, $0.844$ and 1
(top row, second row, third row and bottom, respectively).

The top two panels again show that $W(A)$ is a disk. In the top left (the real case), the three eigenvalues 
of $A$ are (nearly) coincident with two of the roots of $p$, approximating a configuration $(p,A)$ where $A$ 
is a Crabb matrix of order 3 (a $3\times 3$ Jordan block) and $p(z)=(z-\lambda)^{2}$, where
$\lambda$ is the eigenvalue of $A$. Zooming in indeed shows that the three computed eigenvalues of $A$ and
two of the roots of $p$ are nearly coincident, while the third root (not shown) has enormous modulus. However,
in the top right (the complex case), the optimal configuration that is approximated is subtly different. The computed
complex triangular matrix $A$ has one upper triangular entry which is much larger than the others, and only one of
the roots of the computed $p$ is close to the eigenvalues of $A$, while the other two (not shown) have very large modulus.
This indicates that the optimal configuration being approximated is $(p,A)$ where $A$ is block diagonal with a 
$2\times 2$ Jordan block $J$ (the Crabb matrix of order 2) and $p(z)=z-\lambda$, where $\lambda$ is the double eigenvalue of
$J$, with the third eigenvalue of $A$ separated from the others but inside the field of values of $J$. Indeed, zooming in
we find that one root of $p$ and two eigenvalues of $A$ are nearly coincident, with the third eigenvalue separated from them,
though not by much.

In the second row of Figure~\ref{fig:n3m3fovs-subplot}, we see a new configuration:
the three eigenvalues of $A$ are well separated but the three roots of $p$ are (nearly) coincident, and
$\|p\|_{W(A)}$ is (nearly) attained at three points on the boundary. Consequently, as we see in Table~\ref{tab:n3m3},
in the real case $\Zeps$ has two points (one real and one complex) while in the complex case $\Zeps$ contains three complex
points. Again, we see from Table~\ref{tab:n3m3} that $\|d\|$ is small in both cases, indicating approximate nonsmooth stationarity.

In the third row, associated with the locally minimal value $0.844$ that was observed earlier, the fields
of values clearly indicate that the associated matrices are nearly block diagonal, with blocks of order 2 and 1. On the left, the 
eigenvalue corresponding to the $1\times 1$ block lies inside the field of values of the $2\times 2$ block, while on the right, 
it does not; however, in both cases the numerator and denominator of the Crouzeix ratio are determined by the (approximate)
$2\times 2$ block;  the other block is ``inactive''. This is the reason why the locally minimal value associated with the 
two lower panels is the same as that found for the case $n=2$, $m=2$. 
In contrast, the fields of values shown in the second row are clearly \emph{not} associated 
with block diagonal matrices:  thus, the associated locally minimal value is \emph{not} a locally minimal value when $n=2$.  
To emphasize this point, we might say that $0.698$ is a \emph{genuine} locally minimal value for $n=3$.

\begin{figure}
\begin{center}
\begin{tabular}{cc}
\includegraphics[scale=0.45,trim={1.5cm 1.5cm 1.4cm 0.5cm},clip]{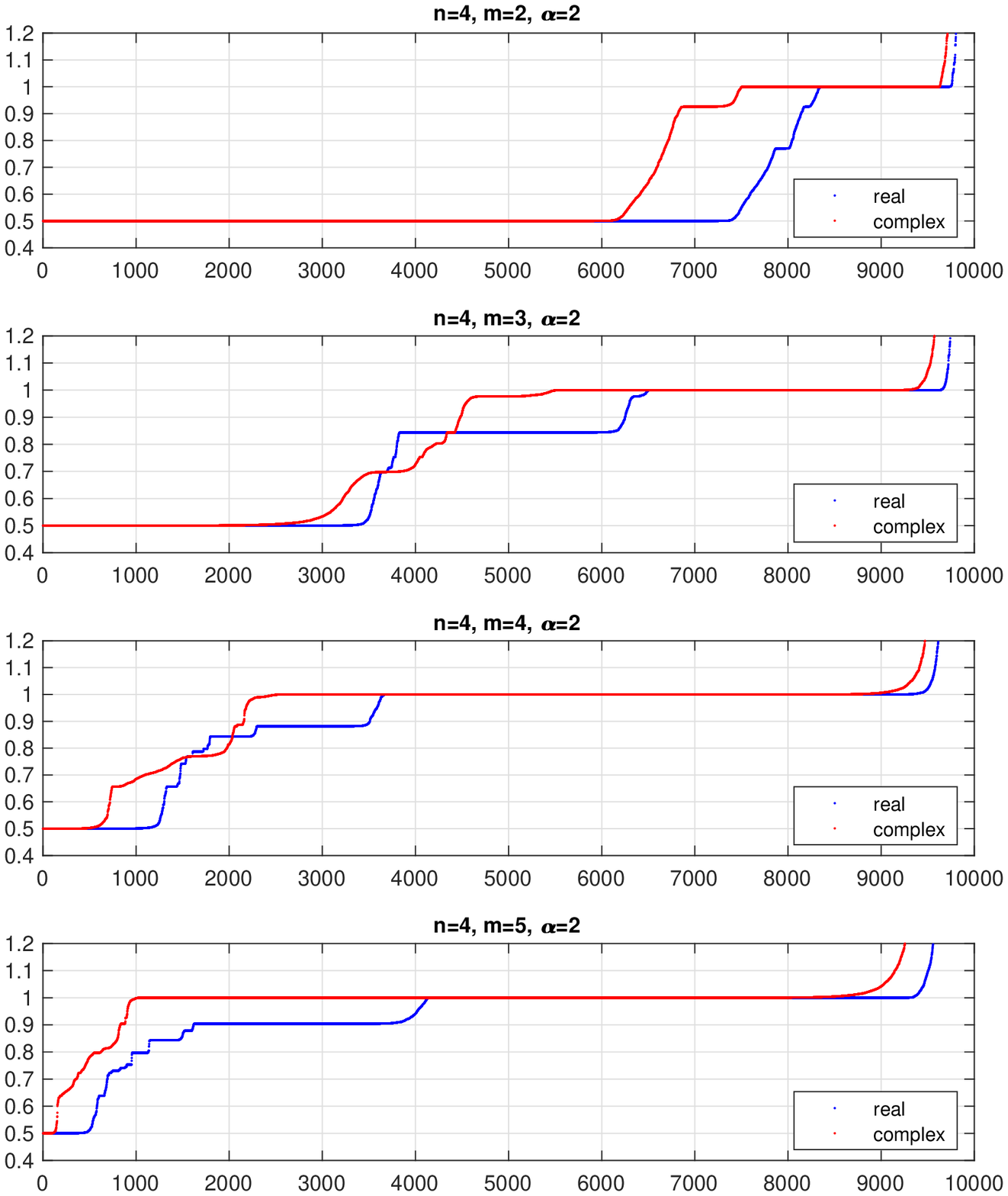} &
\includegraphics[scale=0.45,trim={1.5cm 1.5cm 1.4cm 0.5cm},clip]{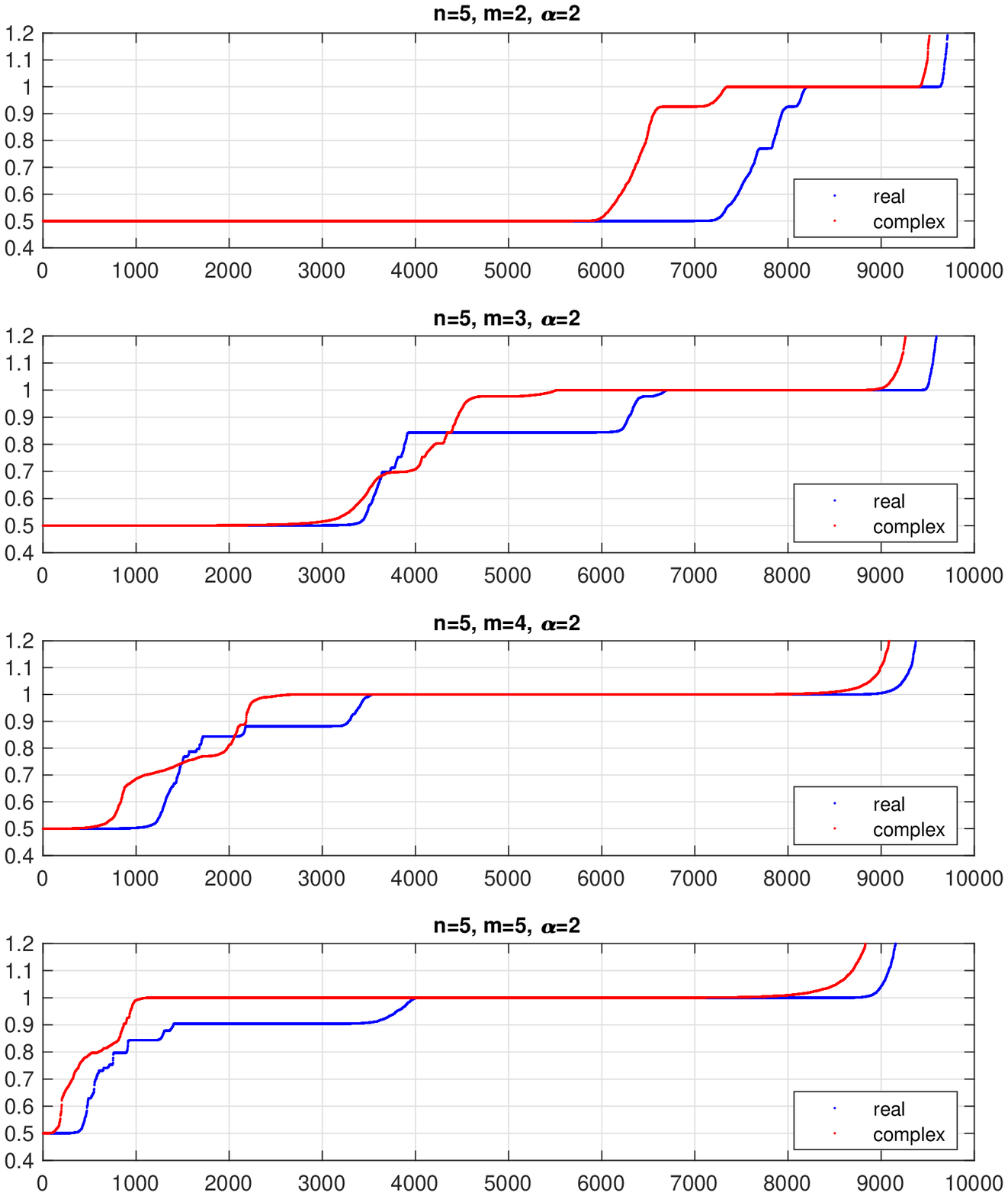}  \\
\end{tabular}
\end{center}
\caption{
Sorted final values of the Crouzeix ratio $f$ obtained for $n=4,5$, $m=2,3,4,5$.}
\label{fig:n45m2345}
\end{figure}

\begin{figure}
\begin{center}
\begin{tabular}{cc}
\includegraphics[scale=0.45,trim={1.5cm 1.5cm 1.4cm 0.5cm},clip]{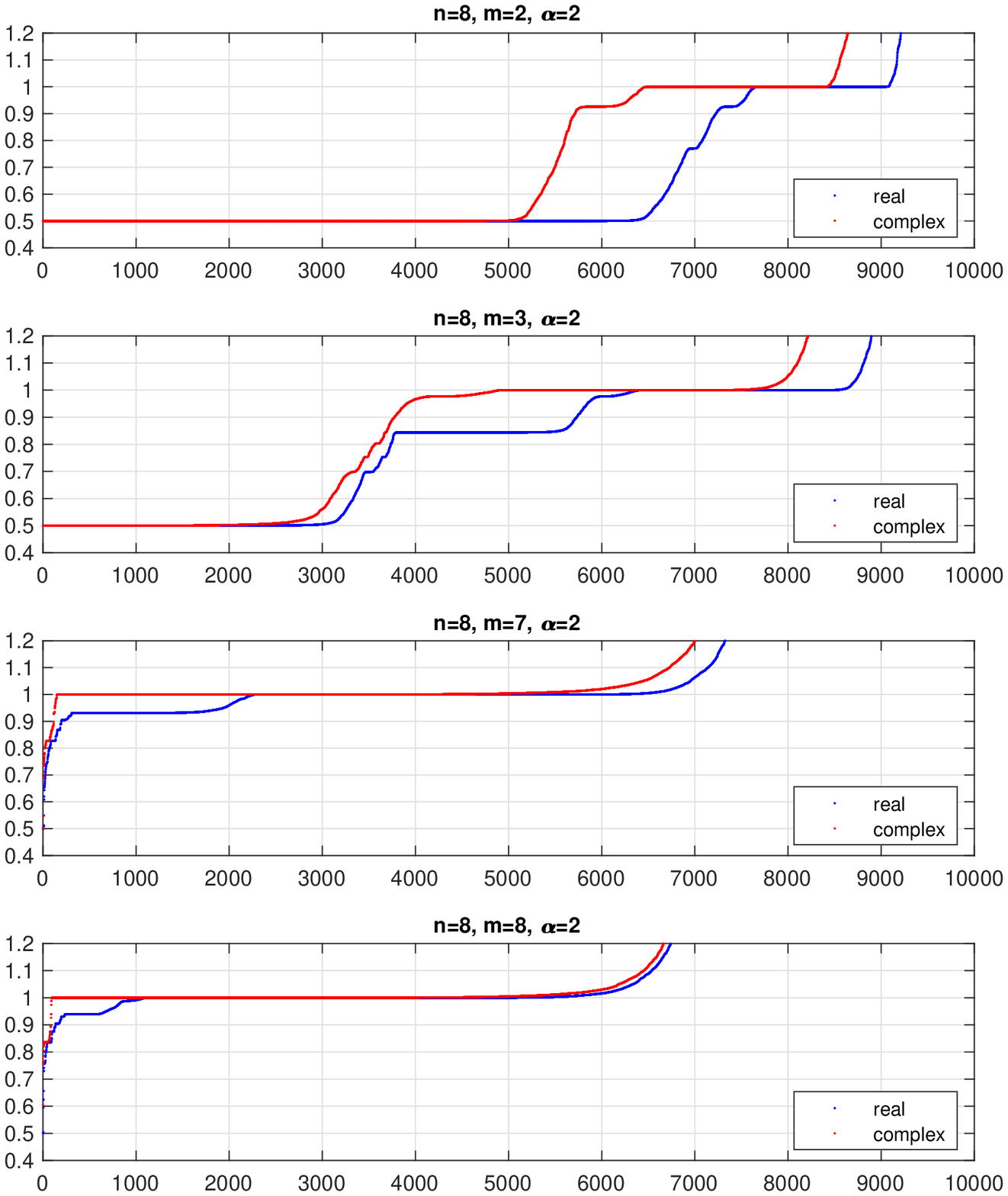} &
\includegraphics[scale=0.45,trim={1.5cm 1.5cm 1.4cm 0.5cm},clip]{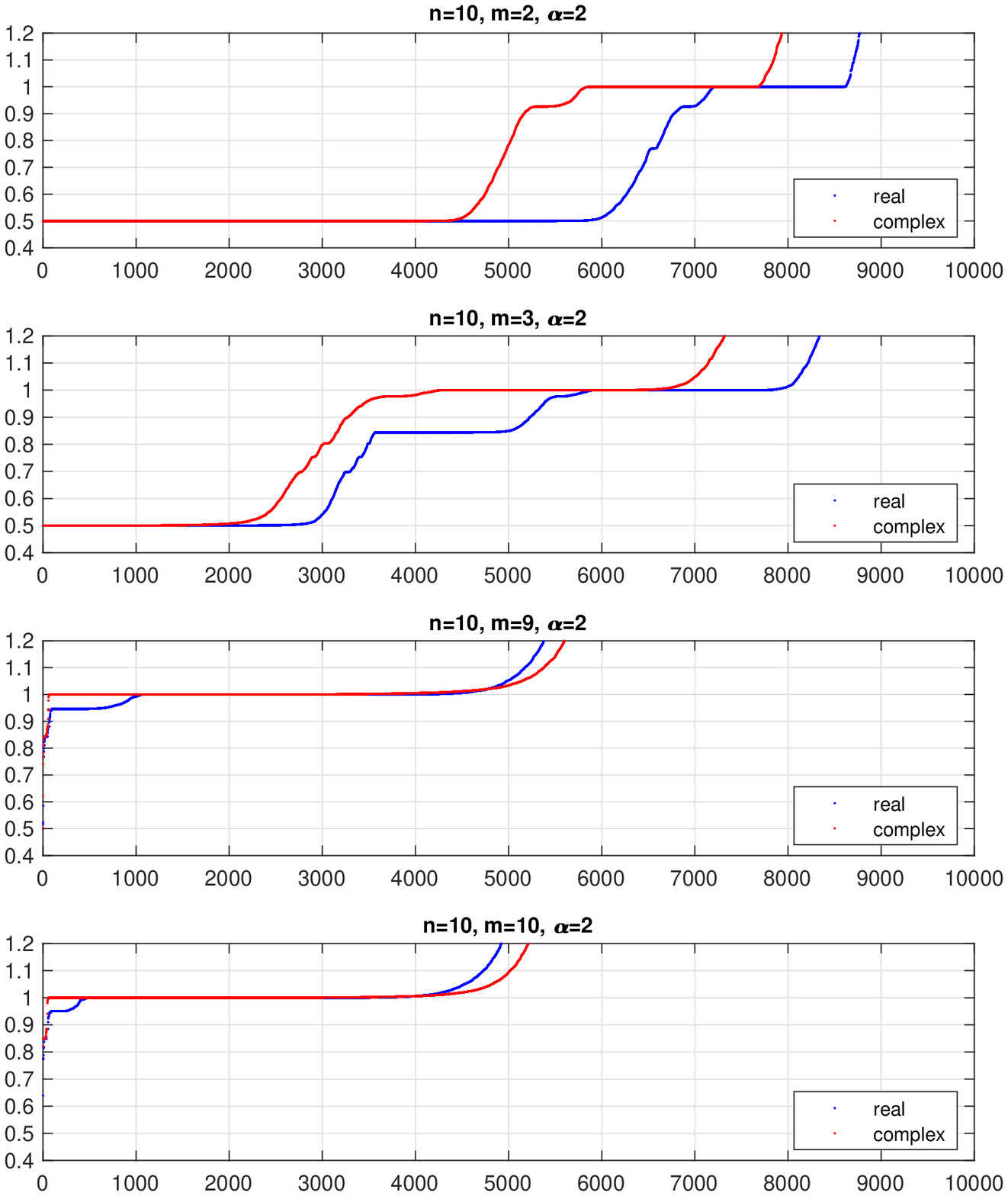}  \\
\end{tabular}
\end{center}
\caption{
Sorted final values of the Crouzeix ratio $f$ obtained for (left) $n=8$, $m=2,3,7$ and 8, and
(right) n=10, $m=2,3,9$ and 10.
}
\label{fig:n8and10}
\end{figure}

Finally, in the fourth row we see two very elongated ice-cream-cone configurations, which are smooth stationary points.

\subsection{Cases with $n\geq 4$}

Figure~\ref{fig:n45m2345} shows results for $n=4$ (left) and $n=5$ (right), with $m=2,3,4,5$. Again we see strong evidence of
locally minimal values between 0.5 and 1, with at least some values coinciding for the real and complex cases. 
Figure~\ref{fig:n8and10} shows results for $n=8$, $m=2,3,7,8$ (left), and n=10, $m=2,3,9,10$ (right). 
For these values of $n$, we can still observe locally minimal values between 0.5 and 1 when $m=2$ or $3$, but for larger $m$
only a handful of results with $f$ below 1 are observed.

\section{Concluding Remarks}\label{sec:conclude}

Crouzeix's conjecture \cite{Cro04} states that the globally minimal value of the Crouzeix ratio~\eqref{ratio} is~0.5, 
regardless of $n$ and $m$, and it was demonstrated in \cite{GreOve18} that 1 is a frequently occurring locally minimal value.
Making use of the heavy-tailed distribution \eqref{heavy} to initialize our optimization runs, 
we have demonstrated for the first time that the Crouzeix ratio has many other locally minimal values between 0.5 and 1,
even for $n=2$, $m=3$ and $n=3$, $m=3$, cases that we studied in detail. Not only did we show that the same function
values are repeatedly obtained for many different starting points, but we also verified that approximate nonsmooth
stationarity conditions hold at computed candidate local minimizers. We also found that the same locally minimal values are often obtained 
both when optimizing over real matrices and polynomials, and over complex matrices and polynomials.
The appearance of so many nonsmooth local minimizers suggests how very complex a function the Crouzeix
ratio is and could perhaps shed some light on why mathematically establishing its global minimum is so difficult.

We think that minimization of the Crouzeix ratio makes a very interesting nonsmooth optimization case study illustrating
among other things how effective the BFGS method is for nonsmooth optimization.
Our method for verifying approximate nonsmooth stationarity is based on what may be a novel approach to finding approximate
subgradients of max functions on an interval, exploiting Chebfun's ability to efficiently find local maximizers on intervals.

Our extensive computations strongly support Crouzeix's conjecture. We have presented results for nearly half a million
optimization runs
reported in Figures \ref{fig:n2m2345}, \ref{fig:n3m2345}, \ref{fig:n45m2345} and \ref{fig:n8and10}, computing the Crouzeix ratio
for about 250 million pairs $(p,A)$.
The computations were done using the high performance computing cluster at New York University, running the code on
hundreds of CPU cores using a {\tt parfor} ``parallel for'' loop in {\sc matlab}. Doing this, even the 80,000 runs
for $n=10$ took less than 3 hours. We always found that the smallest locally minimal value was 0.5.

\bibliography{refs}

\begin{thebibliography}{DHT14}

\bibitem[AO21]{AslOve21}
A.~Asl and M.L. Overton.
\newblock Behavior of limited memory {BFGS} when applied to nonsmooth functions
  and their {N}esterov smoothings.
\newblock In M.~Al-Baali, L.~Grandinetti, and A.~Purnama, editors, {\em Recent
  Developments in Numerical Analysis and Optimization}. Springer, 2021.

\bibitem[BLO02]{BurLewOveMOR}
J.V. Burke, A.S. Lewis, and M.L. Overton.
\newblock Approximating subdifferentials by random sampling of gradients.
\newblock {\em Math.\ Oper.\ Res.}, 27:567--584, 2002.

\bibitem[BLO05]{BurLewOveGradSamp}
J.V. Burke, A.S. Lewis, and M.L. Overton.
\newblock A robust gradient sampling algorithm for nonsmooth, nonconvex
  optimization.
\newblock {\em SIAM Journal on Optimization}, 15:751--779, 2005.

\bibitem[Cla75]{Cla75}
Frank~H. Clarke.
\newblock Generalized gradients and applications.
\newblock {\em Trans. Amer. Math. Soc.}, 205:247--262, 1975.

\bibitem[CP17]{CroPal17}
M.~Crouzeix and C.~Palencia.
\newblock The numerical range is a {$(1+\sqrt{2})$}-spectral set.
\newblock {\em SIAM J. Matrix Anal. Appl.}, 38(2):649--655, 2017.

\bibitem[Cra71]{Cra71}
Michael~J. Crabb.
\newblock The powers of an operator of numerical radius one.
\newblock {\em Michigan Math. J.}, 18:253--256, 1971.

\bibitem[Cro04]{Cro04}
M.~Crouzeix.
\newblock Bounds for analytical functions of matrices.
\newblock {\em Integral Equations and Operator Theory}, 48:461--477, 2004.

\bibitem[Cro16]{Cro16}
Michel Crouzeix.
\newblock Some constants related to numerical ranges.
\newblock {\em SIAM J. Matrix Anal. Appl.}, 37(1):420--442, 2016.

\bibitem[DHT14]{Chebfun14}
T.~A. Driscoll, N.~Hale, and L.~N. Trefethen.
\newblock {\em Chebfun Guide}.
\newblock Pafnuty Publications, Oxford, 2014.

\bibitem[GLO17]{GreLewOve17}
Anne Greenbaum, Adrian~S. Lewis, and Michael~L. Overton.
\newblock Variational analysis of the {C}rouzeix ratio.
\newblock {\em Math. Program.}, 164(1-2, Ser. A):229--243, 2017.

\bibitem[GO18]{GreOve18}
Anne Greenbaum and Michael~L. Overton.
\newblock Numerical investigation of {C}rouzeix's conjecture.
\newblock {\em Linear Algebra Appl.}, 542:225--245, 2018.

\bibitem[Gol77]{Gol77}
A.A. Goldstein.
\newblock Optimization of {L}ipschitz continuous functions.
\newblock {\em Mathematical Programming}, 13:14--22, 1977.

\bibitem[GOS15]{GugOveSte15}
N.~Guglielmi, M.L. Overton, and G.~W. Stewart.
\newblock An efficient algorithm for computing the generalized null space
  decomposition.
\newblock {\em SIAM J. Matrix Anal. Appl.}, 36(1):38--54, 2015.

\bibitem[HJ91]{HorJoh91}
R.A. Horn and C.R. Johnson.
\newblock {\em Topics in Matrix Analysis}.
\newblock Cambridge University Press, Cambridge, U.K., 1991.

\bibitem[Kip51]{Kip51}
R.~Kippenhahn.
\newblock \"{U}ber den {W}ertevorrat einer {M}atrix.
\newblock {\em Math. Nachr.}, 6:193--228, 1951.
\newblock English translation by P.F. Zachlin and M.E. Hochstenbach, Linear and
  Multilinear Algebra 56, pp.\ 185-225, 2008.

\bibitem[LO13]{LewOveNSOquasi}
A.S. Lewis and M.L. Overton.
\newblock Nonsmooth optimization via quasi-{N}ewton methods.
\newblock {\em Math. Program.}, 141(1-2, Ser. A):135--163, 2013.

\bibitem[RW98]{RocWet98}
R.~Tyrrell Rockafellar and Roger J.-B. Wets.
\newblock {\em Variational Analysis}, volume 317 of {\em Grundlehren der
  Mathematischen Wissenschaften [Fundamental Principles of Mathematical
  Sciences]}.
\newblock Springer-Verlag, Berlin, 1998.

\end{thebibliography}
\bibliographystyle{alpha}

\end{document}